\documentclass[12pt]{article}
\usepackage{amsmath,amsfonts,amsthm,amssymb}
\usepackage{latexsym,amsmath}
\usepackage{graphicx,psfrag}

\numberwithin{equation}{section}

 \DeclareMathOperator{\id}{id}

 \DeclareMathOperator{\Fix}{Fix}
 \DeclareMathOperator{\intr}{int}

 \newcommand{\R}{\mathbb{R}}

 \newtheorem{thm}{Theorem}[section]
 \newtheorem{lem}[thm]{Lemma}
 \newtheorem{prop}[thm]{Proposition}

 \theoremstyle{definition}
 \newtheorem{defn}[thm]{Definition}
 \newtheorem{exam}[thm]{Example}
 \newtheorem{rem}[thm]{Remark}
 \theoremstyle{remark}

\begin{document}

\title
{On tame embeddings of solenoids \\ into 3-space} \maketitle
\begin{center}
BOJU JIANG, jiangbj@math.pku.edu.cn
\\
{\it Depart. of Mathematics, Peking University, Beijing  100871,
China\\}

\end{center}

\begin{center}
SHICHENG WANG, wangsc@math.pku.edu.cn\\
{\it Depart. of Mathematics, Peking University, Beijing 100871,
China\\}
\end{center}

\begin{center}
HAO ZHENG, zhenghao@mail.sysu.edu.cn\\
{\it Depart. of Mathematics, Zhongshan University, Guangzhou
510275, China\\}
\end{center}

\begin{center}
QING ZHOU, qingzhou@sjtu.edu.cn\\
{\it Depart. of Mathematics, Jiaotong University, Shanghai 200030,
China\\}
\end{center}

\begin{abstract}

Solenoids are ``inverse limits'' of the circle, and the classical
knot theory is the theory of tame embeddings of the circle into the
3-space. We give some general study, including certain
classification results, of tame embeddings of solenoids into the
3-space as the ``inverse limits'' of the tame embeddings of the
circle.

Some applications are discussed. In particular, there are ``tamely''
embedded solenoids $\Sigma\subset \R^3$ which are strictly achiral.
Since solenoids are non-planar, this contrasts sharply with the
known fact that if there is a strictly achiral embedding $Y\subset
\R^3$ of a compact polyhedron $Y$, then $Y$ must be planar.

\end{abstract}

\maketitle

\section{Introduction and motivations}\label{s:1}

The classical knot theory is the theory of tame embeddings of the
circle into the 3-space, which has become a central topic in
mathematics. The classical theory of knots has many generalizations
and variations: from the circle to graphs, from the circle to higher
dimensional spheres, from tame embeddings to wild embeddings, and so
no. In the present note, we try to setup a beginning of another
generalization: the tame embeddings of solenoids into the 3-space.
In such a study topology and dynamics interact well.

The solenoids may be first defined in topology by Vietoris in 1927
for 2-adic case [V] and by many others later for general cases,
and introduced into dynamics by Smale in 1967 [S]. Solenoids can
be presented either in a rather geometric way (intersections of
nested solid tori, see Definition \ref{defn:nested intersection})
or in a rather algebraic way (inverse limits of self-coverings of
the circle, see Definition \ref{defn:inverse limit}), or in rather
dynamics way (mapping tori over the Cantor set, see [Mc]).

The precise definition of tame embedding of solenoids into the
3-space $\R^3$ will be given in \S 2, but the intuition is quite
naive. Recall we identify $S^1$ with the centerline of the solid
tours $S^1\times D^2$, and say an embedding $S^1\subset\R^3$ is tame
if the embedding can be extended to an embedding $S^1\times
D^2\subset\R^3$. Similarly we consider a solenoid $\Sigma$ as the
nested intersections of solid tori (the defining sequence of
$\Sigma$), and say an embedding $\Sigma\subset\R^3$ is tame if the
embedding can be extended to an embedding of those solid tori into
$\R^3$.

\medskip
Solenoids themselves usually are considered as "wild" set. What
motivated us originally to study the tame embeddings of solenoids
was trying to find a non-planar set which admits a strictly
achiral embedding into the 3-space.

An embedding $A \subset \R^3$ is called strictly achiral, if $A$
stays in the fixed point set of an orientation reversing
homeomorphism $r : \R^3 \to \R^3$. Obviously any planar set has a
strictly achiral embedding. Indeed there is a simple relation
between the two notions of achirality and planarity for compact
polyhedra: If there is a strictly achiral embedding $Y\subset \R^3$
of a compact polyhedron $Y$, then $Y$ must be planar [JW]. In the
sequel, it is natural to ask whether the relation from [JW] still
holds when compact polyhedra are replaced by continua (i.e. compact,
connected metric spaces).

The solenoids are promising, and indeed are proved to be in \S 4,
counterexamples to the question, because on the one hand they are
continua realized as inverse limits of planar sets, on the other
hand they are non-planar themselves (see [Bin], and [JWZ] for a
short proof). In order to design a strictly achiral embedding
$\Sigma \subset \R^3$ for a solenoid $\Sigma$, we need careful and
deep discussions about the tame embeddings of solenoids.

Among other motivations, when a solenoid $\Sigma \subset S^3$ is
realized as a hyperbolic attractor of a dynamics (called Smale
solenoid), the embedding $\Sigma \subset S^3$ is automatically tame.
Our study also gives some application on this aspect.

\medskip

The contents of the paper are as follows.

In \S 2, we give the precise definitions of tame solenoids and
related notions. We present a lemma about convergence of
homeomorphisms which will be repeatedly used in the paper. With
this lemma we give an alternative description of the tame
solenoids via the language of mapping torus. For comparison we
also construct examples of non-tame embeddings of solenoids vis
the mapping torus description.

In \S 3, we classify the tame solenoids in the 3-sphere. The
classification is based on an important notion, the ``maximal''
defining sequences of tame solenoids. We also give applications
including (1) the knotting, linking and invariants of tame
solenoids; (2) there are uncountably many unknotted 2-adic tame
solenoids; (3) the number of Smale solenoids $\Sigma\subset S^3$
realizing the periodic solenoids of type $(w_1, ..., w_k)$ is finite
if all $w_n\le 3$, and is countably infinite otherwise.

In \S 4, we give criterions of when a tame solenoid in the
3-sphere is achiral or strictly achiral in term of its defining
sequence,  and construct strictly achiral tame solenoids to fulfil
our original motivation. Indeed we give a simple criterion when a
solenoid has a strictly achiral embedding into the 3-sphere.

\medskip

All terminologies not defined are standard. For 3-manifolds, see
[Ja]; for knot theory see [A]; and for braid theory, see [Bir].

\section{Tame embeddings of solenoids, Preliminaries}\label{s:2}

\subsection{Definitions of solenoids and their tame embeddings}

Let $N=D^2\times S^1$ be the solid torus, where $D^2$ is the unit
disc and $S^1$ is the unit circle. Then $N$ admits a standard
metric. A {\it meridian disk} of $N$ is a $D^2$ slice of $N$. A {\it
framing} of $N$ is a circle on $\partial{N}$ which meets each
meridian disk of $N$ at exactly one point.

\begin{defn}\label{defn:inverse limit}
(1) For a sequence of maps $\{\phi_n: X_n\to X_{n-1}\}_{n\ge1}$
between continua, the {\it inverse limit} is defined to be the
subspace
$$\Sigma = \{(x_0,x_1,...,x_n,...) \mid x_n\in X_n, \; x_{n-1}=\phi_n(x_n)\}$$
of the product space $\Pi_{n=0}^\infty X_n$.

(2) The inverse limit of a sequence of covering maps $\{\phi_n: S^1
\to S^1\}_{n\ge1}$, where $\phi_n$ is of degrees $w_n\ne0$, is
called a {\it solenoid} of type $\varpi=(w_1,w_2,...,w_n,...)$.
\end{defn}

\begin{defn}\label{defn:nested intersection}
(1) Call an embedding  $e:N\to \intr N$, or simply call the image
$e(N)$, a {\it thick braid of winding number $w$} if $e$ preserves
the $D^2$-fiberation and descends to a covering map $S^1 \to S^1$
given by $e^{it} \mapsto e^{iwt}$. Note that the composition of
finitely many thick braids is also a thick braid.

(2) Let $\{e_n:N\to N\}_{n\ge1}$ be an infinite sequence of thick
braids of winding numbers $w_n\ne0$. Let $\psi_n = e_n \circ \cdots
\circ e_2 \circ e_1$ and $N_n=\psi_n(N)$. Then we have an infinite
sequence $N = N_0 \supset N_1 \supset N_2 \supset \cdots \supset N_n
\supset \cdots$ of thick braids. If the diameters of the meridian
disks of $N_n$ tend to zero uniformly as $n\to\infty$ then we call
$\Sigma = \cap_{n\ge0}\psi_n(N) = \cap_{n>0}N_n$ a {\it solenoid} of
type $\varpi=(w_1,w_2,...,w_n,...)$.
\end{defn}

There is  quite a rich theory about solenoids developed in
1960-1990's. We just list some basic facts (see [Mc], [R] and
references therein) as

\begin{thm}
(1) The above two definitions of solenoids are equivalent; each
solenoid is determined by its type $\varpi$.

(2) Two solenoids $\Sigma$ and $\Sigma'$ of types $\varpi$ and
$\varpi'$ respectively are homeomorphic if deleting finitely many
terms from $\varpi$ and $\varpi'$ can make them identical. Moreover,
$\Sigma$ is the circle if and only if all except finitely many $w_n$
are equal to $\pm1$.

(3) Each solenoid $\Sigma$ is connected, compact and has topological
dimension one. Moreover, if $\Sigma$ is not the circle, then
$\Sigma$ has uncountably many path components.
\end{thm}

We assume below that all winding numbers involved in the definition
of solenoid are greater than 1, unless otherwise specified. In
particular, a solenoid is not the circle.

\begin{defn}\label{defn:tameness}
(1) In Definition \ref{defn:nested intersection}, call
$\{N_n\}_{n\ge0}$ a {\it defining sequence} of the solenoid
$\Sigma$, and call $\Sigma\subset N$ a {\it standard embedding} of
$\Sigma$ in the solid torus $N$.

(2) $\Sigma\subset N$ is called a {\it tame} embedding of $\Sigma$
in the solid torus $N$, if there is a homeomorphism $f: (N,\Sigma
)\to (N,\Sigma')$ for some standard embedding $\Sigma'\subset N$;
then call $\{f^{-1}(N_n')\}_{n\ge0}$ a {\it defining sequence} of
$\Sigma$, where $\{N_n'\}_{n>0}$ is a defining sequence of
$\Sigma'$.

(3) An embedding $\Sigma \subset S^3$ of solenoid is called $\it
tame$, if the embedding can be factored as $\Sigma \subset N
\subset S^3$ in which $\Sigma \subset N$ is tame; then each
defining sequence $\{N_n\}_{n\ge0}$ of $\Sigma\subset N$ is also
considered as a defining sequence of $\Sigma\subset S^3$, and we
have $S^3 \supset N=N_0 \supset N_1 \supset N_2 \supset \cdots
\supset N_n \supset \cdots \supset \Sigma$.
\end{defn}

\begin{rem}\label{coherent}
From Definition \ref{defn:tameness}, for  each  defining sequence
$\{N_n\}_{n\ge0}$ of a tame solenoids in $S^3$, we always assume
that the $D^2$-slices of all $N_i$ are coherent.
\end{rem}

\begin{defn}\label{defn:equiv}
Call two tame solenoids $\Sigma, \Sigma'\subset S^3$ {\it
equivalent} if there is an orientation preserving homeomorphism $f:
S^3 \to S^3$ such that $f(\Sigma) = \Sigma'$.
\end{defn}

\begin{defn}\label{defn:equivseq}
Say two defining sequences $\{N_n\}_{n\ge0}$ and $\{N'_n\}_{n\ge0}$
of tame solenoids in $S^3$ are {\it strongly equivalent}, if there
are orientation preserving homeomorphism $f_0 : (S^3,N_0) \to
(S^3,N'_0)$ and orientation preserving homeomorphisms $f_n :
(N_{n-1},N_n) \to (N'_{n-1},N'_n)$ with $f_n|_{\partial{N}_{n-1}} =
f_{n-1}|_{\partial{N}_{n-1}}$ for $n\ge1$. Say $\{N_n\}_{n\ge0}$ and
$\{N'_n\}_{n\ge0}$ are {\it equivalent}, if $\{N_{k+n}\}_{n\ge0}$
and $\{N'_{k'+n}\}_{n\ge0}$ are strongly equivalent for some
non-negative integers $k$ and $k'$.
\end{defn}

\begin{rem}
(1) A defining sequence $\{N_n\}_{n\ge0}$ of a tame solenoid
$\Sigma\subset S^3$ carries the information of the braiding of $N_n$
in $N_{n-1}$ and the knotting of $N_n$ in $S^3$. The winding numbers
$w_n$, the simplest invariant of the braiding of $N_n$ in $N_{n-1}$,
give rise to the type of the abstract solenoid $\Sigma$.

(2) Suppose $\Sigma\subset S^3$ is a tame embedding given by a
defining sequence $\{N_n\}_{n\ge0}$, then any infinite subsequence
of $\{N_n\}_{n\ge0}$ is a defining sequence of the same embedding.

(3) Definitions \ref{defn:tameness}-\ref{defn:equivseq} also apply
to other 3-manifolds in the obvious way.
\end{rem}

In this paper we view each braid $\beta$ is defined in $D^2\times
[0,1]$ with ends stay in $D^2\times \{0,1\}$. The closure
$\overline\beta \subset D^2\times S^1$ is obtained by identifying
$D^2\times 0$ and $D^2\times 1$ via the identity. Conversely, for
each closed braid $\overline\beta \subset D^2\times S^1$, cutting
$D^2\times S^1$ open along a $D^2$ slice yields a braid
$\beta\subset D^2\times [0,1]$ up to conjugacy. Therefore, we have
the 1-1 correspondence between the set of closed braids of winding
number $n$ and the set of conjugacy classes of braids on $n$
strands.

Note that a tubular neighborhood of a closed braid in $D^2\times
S^1$ is a thick braid and, conversely, a framing of a thick braid in
$D^2\times S^1$ is a closed braid.

\subsection{Convergence of homeomorphisms}

The following lemma will be used repeatedly in this paper.

\begin{lem}\label{lem:converge}
Let $X,Y$ be compact metric spaces and $\{ f_n : X \to Y
\}_{n\ge0}$ be a series of homeomorphisms. If there exist subsets
$\{U_n\}_{n\ge0}$ of $X$ and positive numbers
$\{\varepsilon_n\}_{n\ge0}$ such that
  \par (1) $U_n \subset U_{n-1}$ and
  $f_n|_{X\setminus U_{n-1}} = f_{n-1}|_{X\setminus U_{n-1}}$,
  \par (2) $\lim_{n\to\infty} \varepsilon_n = 0$ and
  $d(f_m(x),f_n(x)) \leq \varepsilon_n$,
  $d(f_m^{-1}(y),f_n^{-1}(y)) \leq \varepsilon_n$
  for $x \in X$, $y \in Y$, $m \ge n$,
  \par (3) Both $\cap_{n\ge0}U_n$ and $\cap_{n\ge0}f(U_n)$
  have no interior points,
  \par \noindent
then $f_n$ uniformly converges to a homeomorphism $f :
(X,\cap_{n\ge0}U_n) \to (Y,\cap_{n\ge0}f_n(U_n))$. Moreover, if in
addition,
  \par (4) $X=Y$ and $d(f_n(x),x) \leq \varepsilon_n$ for $x \in U_n$,
  \par\noindent
then $\cap_{n\ge0} U_n$ lies in the fixed point set of $f$.
\end{lem}

\begin{proof}
Since both $X$ and $Y$ are compact, it follows from Condition (2)
that $f_n$ and $f_n^{-1}$ uniformly converge to continuous maps $f
: X \to Y$ and $g : Y \to X$, respectively. By Condition (1), we
have $f|_{X\setminus U_n} = f_n|_{X\setminus U_n}$ and
$g|_{X\setminus f(U_n)} = f_n^{-1}|_{X\setminus f(U_n)}$ hence
$f|_{X\setminus \cap_{n\ge0}U_n}$ is the inverse of
$g|_{Y\setminus \cap_{n\ge0}f_n(U_n)}$.

Since $X\setminus \cap_{n\ge0}U_n$ is dense in $X$ by Condition
(3), on which the continuous map $gf : X \to X$ acts as the
identity, it follows that $gf=\id_X$. Similarly, we have
$fg=\id_Y$. So the conclusion follows.

The ``Moreover'' part is clear.
\end{proof}

\subsection{Tame solenoids as mapping tori}

\begin{defn}\label{defn:nested disks}
Let $C \subset D^2$ be a Cantor set. Say an orientation preserving
homeomorphism $f : (D^2,C) \to (D^2,C)$ {\it factors through nested
disks} if there exist subsets $\{U_n\}_{n\ge0}$ of $D^2$ such that
$U_{n+1} \subset \intr U_n$, $\cap_{n\ge0}U_n=C$ and each $U_n$ is a
disjoint union of 2-disks on which $f$  cyclically permutes the
components of $U_n$.
\end{defn}

\begin{thm}\label{thm:alt_defn}
A solenoid $\Sigma \subset S^3$ is tame if and only if $\Sigma$
has a neighborhood $N$ such that $(N,\Sigma)$ is homeomorphic to
the mapping torus of a homeomorphism $f : (D^2,C) \to (D^2,C)$
which factors through nested disks.
\end{thm}

The Proof depends on the following lemma which will be also used in
several other places of the paper.

\begin{proof}[Proof of Theorem \ref{thm:alt_defn}]
Sufficiency is clear: Suppose $f : (D^2,C) \to (D^2,C)$ is an
orientation preserving homeomorphism which factors through nested
disks $\{U_n\}_{n\ge0}$. Then for each $n\ge0$, let $N_n =
U_{n}\times[0,1]/f$. Then $N_{n+1}$ is a thick braid in $N_n$ and
$\cap_{n\ge0}N_n=C\times [0,1]/f \subset N_0$ is a tame embedding
of solenoid with defining sequence $\{N_n\}_{n\ge0}$.

Necessity. Let $\Sigma \subset N$ be a tame embedding with
defining sequence $\{N_n\}_{n\ge0}$. We will show that
$(N,\Sigma)$ is homeomorphic to the mapping torus of a
homeomorphism $f : (D^2,C) \to (D^2,C)$ which factors through
nested disks. Fix a meridian disk $D$ of $N$. For $n\ge0$, let
$U_n=D \cap N_n$. Clearly $U_n$ are disjoint disks,
$U_{n+1}\subset \intr U_n$, and $C=\cap_{n\ge0}U_n$ is a Cantor
set in $D^2$.

Now recursively define homeomorphism $f_n : (D,U_n) \to (D,U_n)$
such that $f_n|_{D\setminus U_{n-1}} = f_{n-1}|_{D\setminus
U_{n-1}}$ and $T_{f_n} \cong (N_0,N_n)$, where $T_{f_n}$ is the
mapping torus of $f_n$. Since $N_n$ is connected, $f_n$ permutes the
components of $U_n$ cyclically. Then recursively define
$D^2$-fiberation preserving homeomorphism $\tilde f_n : T_{f_n} \to
(N_0,N_n)$ such that $\tilde f_n$ is identical to $\tilde{f}_{n-1}$
on the mapping torus over ${D\setminus U_{n-1}}$. Since the
diameters of the components of $U_n$ and the meridian disks of $N_n$
tend to zero uniformly as $n\to\infty$ and $\tilde f_n$ is
$D^2$-fiberation preserving, by Lemma \ref{lem:converge} $f_n$
uniformly converges to a homeomorphism $f : (D,D\cap\Sigma) \to
(D,D\cap\Sigma)$ and $\tilde{f}_n$ uniformly converges to a
homeomorphism $\tilde{f} : T_f \to (N_0,\Sigma)$. It is clear that
$f$ factors through $\{U_n\}_{n\ge0}$.
\end{proof}

\subsection{Non-tame embeddings of solenoids}

As there are non-tame embeddings of the circle into the 3-space,
there are non-tame embeddings of solenoids. In fact, any tame
embedding of solenoid can be modified in a simple way to a non-tame
embedding, which we illustrate below by a concrete example.

Suppose nested disks $\{U_n\}_{n\ge0}$ and an orientation preserving
homeomorphism $f : (D^2,C) \to (D^2,C)$ which  factors through
$\{U_n\}_{n\ge0}$ have been chosen as in previous subsection so that
each $U_n$ consists of $2^n$ disks. Label those disks by $U_1=D_0
\cup D_1$, $U_2= D_{00}\cup D_{01} \cup D_{10} \cup D_{11}$ with
$D_{ij}\subset D_i$, where $i,j \in \{0,1\}$, and so on. For
simplicity, denote by $D'_n,D''_n$ the components
$D_{11...10},D_{11...11}  \subset U_n$ respectively.

\bigskip
\begin{center}
    \psfrag{b}[]{$D'_1 \times [0,1]$}
    \psfrag{c}[]{$D'_2 \times [0,1]$}
    \psfrag{B}[]{$e_1(D'_1 \times [0,1])$}
    \psfrag{C}[]{$e_2(D'_2 \times [0,1])$}
    \psfrag{p}[]{$p$}
    \includegraphics[scale=.8]{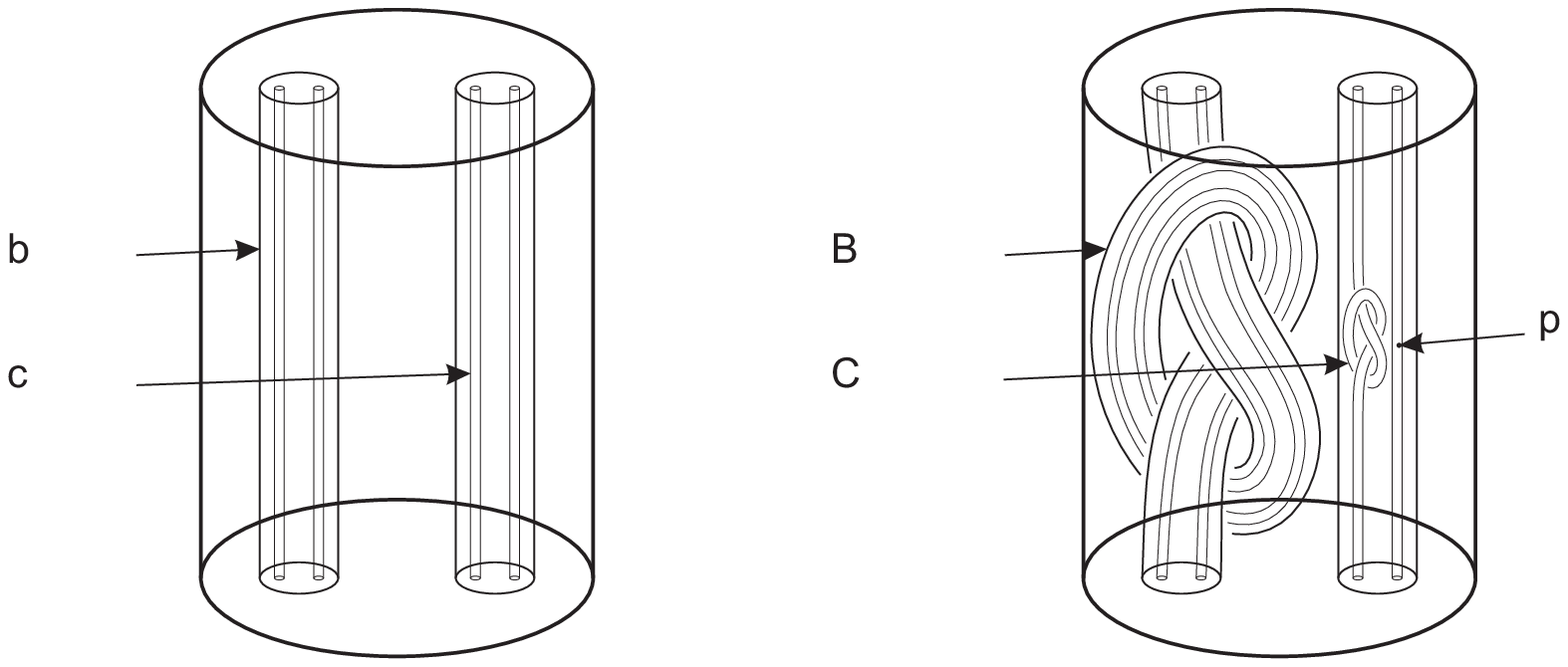}
\end{center}
\centerline{Figure 1}
\bigskip

In the left side of Figure 1, we see the embedding sequence
$D^2\times [0,1]\supset U_1\times [0,1]\supset U_2\times
[0,1]\supset \cdots \supset C\times [0,1]$. We define a re-embedding
$e :C\times [0,1] \subset D^2\times [0,1]$ relative to $D^2\times
\{0,1\}$ by assembling the re-embeddings $\{e_n : D'_n\times
[0,1]\subset D^2\times [0,1]\}_{n\ge1}$, where $e_n$ is constructed
by replacing an unknotted middle portion of $D'_n\times [0,1]$ by a
knotted one such that the knotted portion of $D'_n\times [0,1]$ is
contained in a 3-ball $B_n\subset (D''_{n-1} \setminus D'_n) \times
[0,1]$, the diameter of $B_n$ tending to zero and $B_n$ converging
to a point $p$. The situation is indicated as the right side of
Figure 1. Note that for each small open 3-ball neighborhood $B$ of
$p$ there is an arc component of $B\cap e(C\times [0,1])$ which is
knotted in $B$. Denote this property by $(*)$.

Let $T_f$ be the mapping torus of the homeomorphism $f : (D^2,C) \to
(D^2,C)$.  Then the quotient of $e(C\times [0,1])\subset D^2\times
[0,1]$ in $T_f$ is an embedding of 2-adic solenoid into the solid
torus $T_f$. This embedding is not tame, since a tame embedding do
not have property $(*)$.

\section{Classification and applications}\label{s:3}

\subsection{Maximal defining sequences}

As can be easily seen, a tame solenoid $\Sigma \subset N$ has a
lot of defining sequences and there is no way to choose a
``minimal'' one from them. However, ``maximal'' defining sequence
of a tame solenoid can be defined in the sense that every other
defining sequence is equivalent to a subsequence of it (see
Proposition \ref{prop:maximal}).

Before giving the precise definition we fix a notation. For a
properly embedded surface $S$ (resp. an embedded 3-manifold $P$) in
a 3-manifold $M$, we use $M \setminus S$ (resp. $M \setminus P$) to
denote the resulting manifold obtained by splitting $M$ along $S$
(resp. removing the interior of $P$).

\begin{defn}
Call a defining sequence $\{N_n\}_{n\ge0}$ of a tame solenoid
$\Sigma\subset S^3$ {\it maximal} if $N_n\setminus N_{n+1}$
contains no essential torus for each $n\ge0$.
\end{defn}

\begin{prop}\label{prop:maximal}
Every tame solenoid $\Sigma \subset S^3$ has a maximal defining
sequence. Moreover, if $\{N_n\}_{n\ge0}$ is a maximal defining
sequence of $\Sigma$ then every defining sequence of $\Sigma$ is
equivalent to a subsequence of $\{N_n\}_{n\ge0}$. Hence all
maximal defining sequences of $\Sigma$ are equivalent.
\end{prop}

The following lemma will be repeatedly used in this paper.

\begin{lem}\label{lem:braid}
Suppose $N'$ is a thick braid in $N$ and $\Gamma$ is the
JSJ-decomposition tori of $N\setminus N'$. Then

(1) each component $T$ of $\Gamma$ bounds a solid torus $N^*$ such
that $N^*$ is a thick braid in $N$ and  $N'$ is a thick braid in
$N^*$;

(2) each component of $(N\setminus N')\setminus \Gamma$ contains no
essential torus;

(3) for each solid torus $N''$ with $N' \subset N'' \subset N$,
$\partial{N}''$ is isotopic in $N \setminus N'$ to a component of
$\Gamma \cup \partial{N} \cup \partial{N'}$.
\end{lem}

\begin{proof}
(1) Let $w$ be the winding number of $N'$ in $N$ and let $D$ be a
meridian disk of $N$ which meets $N'$ at $w$ meridian disks of
$N'$. Then $(N\setminus N')\setminus D \cong P_{w}\times I$, where
$P_{w}$ is the $w$-punctured disc. Isotope $T$ so that $T\cap D$
has minimum number of components. Then a standard argument in
3-manifold topology shows that each component of $T \setminus D$
in $P_{w}\times I$ is a vertical annulus which separates a
vertical $D^2\times I$ from $N \setminus D$. Therefore $T$ bounds
a solid torus $N^*$ which is a thick braid in $N$ and clearly $N'$
is a thick braid in $N^*$.

(2) Let $Q$ be a component of $(N\setminus N')\setminus \Gamma$.
By JSJ theory [Ja], $Q$ is either simple hence contains no
essential torus by definition, or a Seifert piece. Suppose $Q$ is
a Seifert piece. Then $Q$ is also a Seifert piece of knot
complement with incompressible boundary. According to [Ja, IX.22.
Lemma], $Q$ is either a torus knot space, or a $P_w \times S^1$
where $P_w$ is the $w$-punctured disc with $w\ge 2$ , or a cable
space. Since $\partial Q$ has at least two components, $Q$ is not
the torus knot space. By the conclusion of (1) and the fact that
$N'$ is connected, one can verify that there is no embedding of
$P_w \times S^1$ in $N \setminus N'$ with incompressible boundary
for $w\ge2$. Therefore $Q$ is a cable space. It is known that a
cable space contains no essential torus.

(3) By the conclusions of (2), it suffices to show that
$\partial{N}''$ is incompressible in $N\setminus N'$. Suppose
$\partial{N}''$ has a compressing disk $D$ in $N\setminus N'$.
Then after surgery on $D$ we will get a separating 2-sphere $S^2$
in $N\setminus N'$ such that each component of $(N\setminus
N')\setminus S^2$ contains a boundary torus, which contradicts
that $N\setminus N'$ is irreducible.
\end{proof}

\begin{proof}[Proof of Proposition \ref{prop:maximal}]
Let $\{N_n\}_{n\ge 0}$ be a defining sequence of $\Sigma \subset
S^3$. Then we have $S^3\supset N=N_0 \supset N_1 \supset N_2
\supset \cdots \supset N_n \supset \cdots \supset
\Sigma=\cap_{n=1}^{\infty} N_n$ and
$$ N_0\setminus N_1 \subset N_0\setminus N_2 \subset \cdots
  \subset N_0\setminus N_n \subset \cdots \subset N_0\setminus\Sigma
 = \cup_{n=1}^{\infty} N_0\setminus N_n.
$$
Moreover, since the winding number of $N_{n+1}$ in $N_n$ is greater
than 1 for every $n$, $\partial{N}_n$ is not parallel to
$\partial{N}_j$ for $n\ne j$.

Now we refine the defining sequence $\{N_n\}_{n\ge 0}$ to a
maximal one. Let $\Gamma_n$ be the JSJ-decomposition tori of $N_n
\setminus N_{n+1}$. Note that the set $\bigcup_{0\le n<j}\Gamma_n
\cup \bigcup_{0<n<j}\partial{N}_n$ is the JSJ-decomposition tori
of $N_0\setminus N_j$. Define $\Gamma(\Sigma) = \bigcup_{n\ge0}
(\Gamma_n\cup\partial{N}_n)$. By Lemma \ref{lem:braid} (1)(2) we
can re-index the components of $\Gamma (\Sigma)$ as
$\{T_n\}_{n\ge0}$ so that
  \par (1) each $T_n$ bounds a solid torus $N_n^*$,
  \par (2) each $N_{n+1}^*$ is a thick braid in $N_{n}^*$,
  \par (3) each $N_n^*\setminus N_{n+1}^*$ contains no essential torus.
  \par \noindent
Clearly $\{N_n\}_{n\ge0}$ is a subsequence of $\{N_n^*\}_{n\ge0}$
and by definition $\{N_n^*\}_{n\ge0}$ is a maximal defining
sequence of $\Sigma$.

Let $\{N'_n\}_{n\ge0}$ be another maximal defining sequence of
$\Sigma\subset S^3$. We shall show that $\{N'_n\}_{n\ge0}$ and
$\{N^*_n\}_{n\ge0}$ are equivalent. Since $\Sigma$ is compact, we
have $N^*_j \subset N'_k \subset N^*_0$ for some big integers
$j,k$. By Lemma \ref{lem:braid} (3) and the construction of
$\{N^*_n\}_{n\ge0}$, $\partial N'_k$ is isotopic in
$N^*_0\setminus N^*_j$ to some $\partial{N}^*_m$. Clearly, the
isotopy automatically sends $N'_k$ to $N^*_m$.

Similarly, one argues that $N'_{k+1}$ can be further isotoped in
$N^*_m$ relative to $\Sigma$ to some $N^*_{m'}$, in which $m'$
must be $m+1$ because $N'_k\setminus N'_{k+1}$ contains no
essential torus, and so on. Hence we verify that
$\{N'_{k+n}\}_{n\ge 0}$ is strongly equivalent to
$\{N^*_{m+n}\}_{n\ge0}$ and the conclusion follows.
\end{proof}

\subsection{Classification of tame solenoids}

\begin{thm}\label{thm:classify}
Let $\Sigma, \Sigma' \subset S^3$ be two tame solenoids. The
following statements are equivalent.

(1) $\Sigma, \Sigma'$ are equivalent.

(2) Some defining sequences of $\Sigma,\Sigma'$ are equivalent.

(3) The maximal defining sequences of $\Sigma,\Sigma'$ are
equivalent.
\end{thm}

\begin{proof}
(2) $\Rightarrow$ (1). Without loss of generality, suppose the
defining sequences $\{N_n\}_{n\ge0}$, $\{N'_n\}_{n\ge0}$ of
$\Sigma,\Sigma'$ are strongly equivalent. By definition there are
orientation preserving homeomorphism $f_0 : (S^3,N_0) \to
(S^3,N'_0)$ and orientation preserving homeomorphisms $f_n :
(N_{n-1},N_n) \to (N'_{n-1},N'_n)$ with $f_n|_{\partial{N}_{n-1}}
= f_{n-1}|_{\partial{N}_{n-1}}$ for $n\ge1$.

By Remark 2.5, we assume the $D^2$-slices of all $N_i$,
(respectively of all $N'_i$), are coherent. Then it is easy to see
that we can first isotopy $f_0 : (S^3,N_0) \to (S^3,N'_0)$ so that
$f_0| : N_0 \to N'_0$ is $D^2$-fiberation preserving, then
inductively to isotopy $f_n : (N_{n-1},N_n) \to (N'_{n-1},N'_n)$
for each $n\ge1$ so that $f_n| : N_n \to N'_n$ is $D^2$-fiberation
preserving, and still $f_n|_{\partial{N}_{n-1}} =
f_{n-1}|_{\partial{N}_{n-1}}$.

To apply Lemma \ref{lem:converge}, we set $U_n=N_n$ and extend $f_n$
onto $S^3$ by setting $f_n|_{S^3\setminus N_{n-1}} =
f_{n-1}|_{S^3\setminus N_{n-1}}$. Clearly Conditions (1) and (3) of
Lemma \ref{lem:converge} are satisfied. Since the diameters of the
meridian disks of $N_n$ and $N'_n$ tend to zero uniformly as
$n\to\infty$ and since $f_n$ is $D^2$-fiberation preserving and
$f_n|_{\partial{N}_{n-1}} = f_{n-1}|_{\partial{N}_{n-1}}$, Condition
(2) of Lemma \ref{lem:converge} is also satisfied. Therefore, by
Lemma \ref{lem:converge} $f_n$ uniformly converges to a
homeomorphism $f : (S^3,\cap_{n\ge0}N_n) \to
(S^3,\cap_{n\ge0}N'_n)$. That is, $\Sigma = \cap_{n\ge0}N_n$ and
$\Sigma' = \cap_{n\ge0}N'_n$ are equivalent.

(3) $\Rightarrow$ (2) is obvious.

(1) $\Rightarrow$ (3). Let $f: S^3 \to S^3$ be an orientation
preserving homeomorphism such that $f(\Sigma) = \Sigma'$. Clearly
for each maximal defining sequence $\{N_n\}_{n\ge0}$ of $\Sigma$,
$\{f(N_n)\}_{n\ge0}$ is a maximal defining sequence of $\Sigma'$
and is equivalent to $\{N_n\}_{n\ge0}$. By Proposition
\ref{prop:maximal}, $\{f(N_n)\}_{n\ge0}$, hence $\{N_n\}_{n\ge0}$,
is equivalent to every maximal defining sequence of $\Sigma'$.
\end{proof}

\subsection{Knotting, linking and invariants}

Thanks to the classification theorem, we can talk about the
knotting, linking and invariants of tame solenoids.

\begin{defn}
Call a tame embedding of solenoid $\Sigma\subset S^3$ with defining
sequence $\{N_n\}_{n\ge0}$ is {\it knotted}, if some defining solid
torus $N_n\subset S^3$ is knotted; otherwise call it {\it
unknotted}.
\end{defn}

Note that for a defining sequence $\{N_n\}_{n\ge0}$ of a tame
solenoid, if $N_n$ is knotted then so is $N_{n'}$ for $n'>n$. It
follows from Theorem \ref{thm:classify} that the notion of knotting
is well defined for an equivalent class of tame solenoids.

\begin{defn}\label{defn:linking}
Let $\Sigma, \Sigma'\subset S^3$ be disjoint tame solenoids with
disjoint defining sequences $\{N_n\}_{n\ge0}$ and $\{N'_j\}_{j\ge
0}$ respectively.

(1) Call $\Sigma, \Sigma'$ {\it algebraically linked} if some
linking number $lk(N_n, N'_j)$ (i.e. the linking number of their
centerlines) is non-zero.

(2) Call $\Sigma, \Sigma'$ {\it linked} if some defining solid tori
$N_n, N'_j$ are linked.
\end{defn}

Since two disjoint tame solenoids $\Sigma, \Sigma'\subset S^3$
always have disjoint defining sequences and since $lk (N_n,
N'_j)\ne0$ implies $lk(N_{n'}, N'_{j'})\ne0$ for all $n'\ge n, j'\ge
j$, by Theorem \ref{thm:classify} again the notion of algebraic
linking is well defined for equivalent classes of tame solenoids.

Similarly the notion of linking is well defined, too. In particular,
$\Sigma$ and $\Sigma'$ are linked if and only if there are no
disjoint 3-balls $B$ and $B'$ such that $\Sigma\subset B,
\Sigma'\subset B'$.

\medskip

To define invariants of tame solenoids, the proposition below will
be of help.

\begin{prop}
Up to strong equivalence, each knotted tame solenoid $\Sigma \subset
S^3$ has a unique maximal defining sequence $\{N_n\}_{n\ge 0}$ such
that $N_0$ is knotted and any other defining sequence
$\{N'_n\}_{n\ge 0}$ with knotted $N'_0$ is a subsequence of
$\{N_n\}_{n\ge 0}$.
\end{prop}

\begin{proof}
Follow the proof of Proposition \ref{prop:maximal}. Assume $N_0$ is
knotted and append further into $\Gamma(\Sigma)$ all such
JSJ-decomposition torus of $S^3\setminus N_0$ that bounds a solid
tours in $S^3$, in which $N_0$ is a thick braid. It is a routine
matter to verify that the resulting maximal defining sequence is
exactly what we want.
\end{proof}

For any knot invariant $I$ (for example, the genus, the Gromov
volume, the Alexander polynomial or the Jones polynomial) one has an
invariant $I$ of tame solenoids as below. For a knotted tame
solenoid $\Sigma\subset S^3$, let $\{N_n\}_{n\ge0}$ be the unique
maximal defining sequence from the above proposition. Then the
infinite sequence $I(\Sigma) = \{ I(N_0), I(N_1) ,..., I(N_n),
...\}$ depends only on the equivalence class of $\Sigma$. If a tame
solenoid $\Sigma \subset S^3$ is unknotted, then for any defining
sequence $\{N_n\}_{n\ge 0}$ of $\Sigma$ the sequence $I(\Sigma) = \{
I(N_0), I(N_1) ,..., I(N_n), ...\}$ is identically trivial, say
$\{0,0,...,0,...\}$, if $I$ is either the genus or the Gromov
volume, or $\{1,1,...,1,...\}$, if $I$ is either the Alexander
polynomial or the Jones polynomial.

In general for given numerical function $g$ and knot invariant $I$,
one may organize the sequence $I(\Sigma)$ into a formal series
$I(\Sigma, g)=\sum_{n=0}^\infty g(n) I(N_n) t^n$. We wonder if
$I(\Sigma, g)$ would have interesting properties for certain $g$ and
$I$ as well as for suitable classes of solenoids.

\subsection{Unknotted 2-adic tame solenoids}

Given an unknotted solid torus $N$ in $S^3$, there are exactly two
kinds of thick braids of winding number two in $N$ that are
unknotted in $S^3$ as shown in Figure 2, where the left one is the
left-handed embedding, which we denote by $-1$, and the right one is
the right-handed embedding, which we denote by $+1$. Then any
maximal defining sequence of a unknotted 2-adic tame solenoids in
$S^3$ can be presented as an infinite sequence of $\pm1$.

Let $Z_2$ be the set of infinite sequences $(a_1,a_2,...,a_n,...)$
of $\pm1$. Two such sequences are said to be {\it equivalent}, if
they can be made identical by deleting finitely many terms. By
Theorem \ref{thm:classify} the equivalence classes of unknotted
2-adic tame solenoids are in 1-1 correspondence to the equivalence
classes of $Z_2$. In particular, there are uncountably many
equivalence classes of unknotted 2-adic tame solenoids.

\bigskip
\begin{center}
    \includegraphics[scale=1]{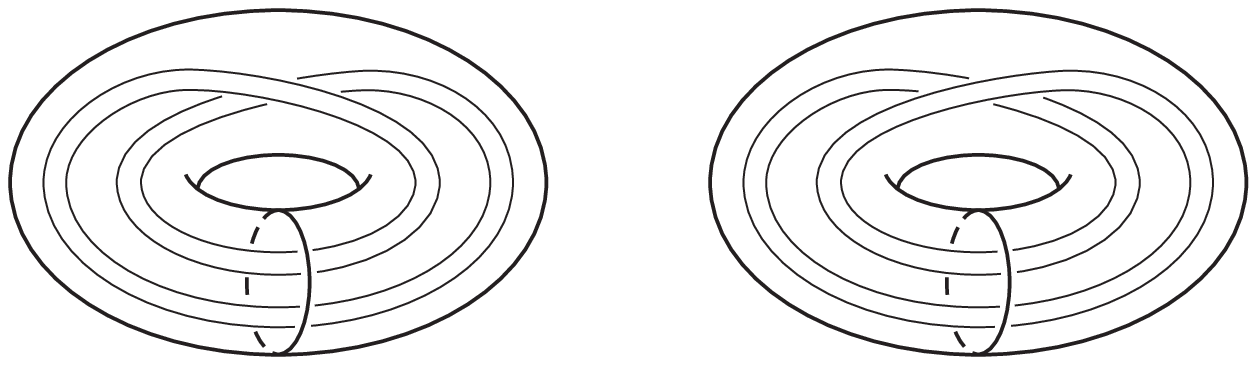}
\end{center}
\centerline{Figure 2}
\bigskip

\subsection{Smale solenoids}

The solenoids were introduced into dynamics by Smale as hyperbolic
attractors in \cite{S}.

\begin{defn}
Let $M$ be a 3-manifold and $f: M\to M$ be a homeomorphism. If there
is a solid torus $N\subset M$ such that $f|_N$ (resp.\ $f^{-1}|_N$)
defines a thick braid as in Definition \ref{defn:nested
intersection} (1), we call the hyperbolic attractor
$\Sigma=\cap_{n=1}^{\infty} f^{n}(N)$ (resp. the hyperbolic repeller
$\Sigma=\cap_{n=1}^{\infty} f^{-n}(N)$) a {\it Smale solenoid}.
\end{defn}

Clearly each Smale solenoid $\Sigma\subset S^3$ is tame. It is known
that a Smale solenoid $\Sigma\subset S^3$ must be unknotted
\cite{JNW}. Moreover, it is proved in \cite{JNW} that if the
non-wondering set $w(f)$ of a dynamics $f$ consists of finitely many
disjoint Smale solenoids, then $w(f)$ consists of two solenoids
(indeed they are algebraically linked).

\begin{defn}
Let $w_1,...,w_k$ be integers greater than $1$. Call a Smale
solenoid $\Sigma\subset S^3$ is of {\it type $(w_1, ..., w_k)$}, if
(1) there is a dynamics $f$ taking $\Sigma$ as an attractor, (2)
there is a defining sequence $\{N_n\}_{n\ge0}$ of $\Sigma$ such that
$f$ sends $N_n$ to $N_{k+n}$ for all $n\ge0$ and (3) $w_n$ is the
winding number of $N_n$ in $N_{n-1}$ for $1\le n\le k$.
\end{defn}

\begin{prop}
Any given type $(w_1, ..., w_k)$ is realized by a Smale solenoid
$\Sigma\subset S^3$. Moreover, the number of Smale solenoids
$\Sigma\subset S^3$ of type $(w_1, ..., w_k)$ is finite if all
$w_n\le 3$, and is countably infinite otherwise.
\end{prop}

\begin{proof}
First, we extend the sequence $(w_1,...,w_k)$ to a infinite one by
setting $w_{k+n}=w_n$. Then choosing an unknotted solid tours
$N_0\subset S^3$ and letting $N_n\subset N_{n-1}$ be a tubular
neighborhood of $\overline{\sigma_1\cdots\sigma_{w_n-1}}$ give rise
to a defining sequence of the desired Smale solenoid, where
$\sigma_i$'s are standard generators of the braid groups.

The ``Moreover'' part follows from the lemma below and Theorem
\ref{thm:classify}.
\end{proof}

\begin{lem}[\cite{MP}]\label{lem:MP}
Let $W_n$ be the set of $n$-strand braids whose closures are
unknotted in $S^3$. Then

(1) $W_n$ has two conjugacy classes as pictured in Figure 2 for
$n=2$;

(2) $W_n$ has three conjugacy classes as pictured in Figure 3 for
$n=3$;

(3) $W_n$ has infinitely many conjugacy classes for $n>3$.
\end{lem}

\bigskip
\begin{center}
    \includegraphics[scale=.8]{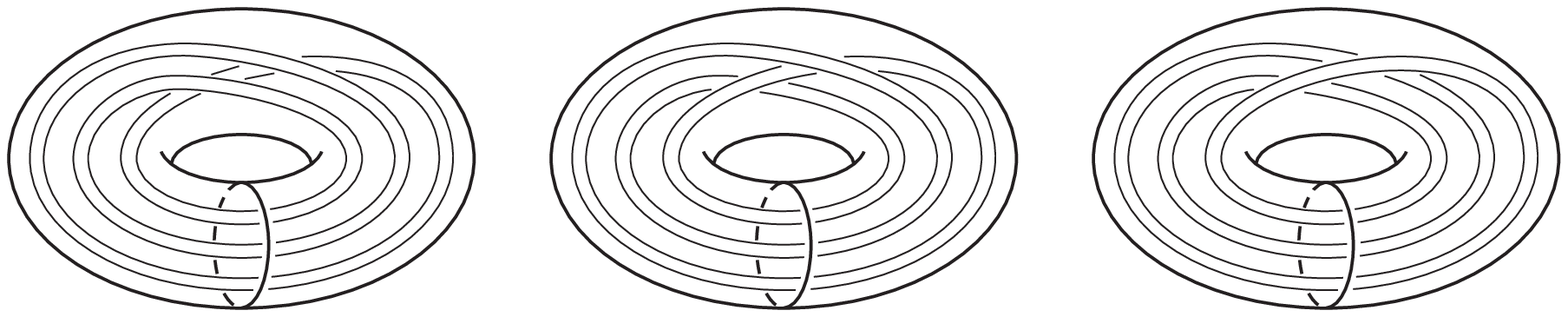}
\end{center}
\centerline{Figure 3}
\bigskip

\section{Chirality of tame solenoids}\label{s:4}

\begin{defn}\label{defn:achiral}
Call a subset $A \subset S^3$ {\it achiral}, if there is an
orientation reversing homeomorphism $r : S^3 \to S^3$ such that
$r(A)=A$. Call $A$ {\it strictly achiral}, if there is an
orientation reversing homeomorphism $r : S^3 \to S^3$ such that
$r(x)=x$ for every $x\in A$.
\end{defn}

In Definition \ref{defn:achiral}, ``achiral'' means setwise achiral,
and ``strictly achiral'' means pointwise achiral. They are two
oppiste extremes among various shades  of achirality in the real
word.

\subsection{Criterions}

By definition, a tame solenoid in $S^3$ is achiral if and only if it
is equivalent to its mirror image. Therefore, by Theorem
\ref{thm:classify} we have the following criterion of the chirality
of tame solenoids.

\begin{thm}\label{thm:achiral}
A tame solenoid given by the maximal defining sequence
$\{N_n\}_{n\ge0}$ is achiral if and only if $\{N_n\}_{n\ge0}$ is
equivalent to its mirror image.
\end{thm}

\begin{exam}
Recall the example in \S 3.4. The mirror image of a maximal defining
sequence of a unknotted 2-adic tame solenoid presented by
$(a_1,a_2,...,a_n,...)$ is presented by $(-a_1,-a_2,...,-a_n,...)$.
Therefore, by Theorem \ref{thm:achiral}, the unknotted 2-adic tame
solenoid presented by $(+1,-1,+1,-1,...)$ is achiral but those
solenoids presented by $(+1,+1,+1,+1,...)$ or $(-1,-1,-1,-1,...)$
are not achiral.
\end{exam}

Below we focus on the strict achirality of tame solenoids. For a map
$f: X\to X$, we use $\Fix(f)$ to denote the fixed point set of $f$.

\begin{defn}
Suppose $A$ is a subset of the solid torus $N$ and $l$ is a given
framing of $N$. Call $A$ is {\it strictly achiral with respect to
$l$} if there exists an orientation reversing homeomorphism $f : N
\to N$ such that $A \cup l \subset \Fix(f)$.
\end{defn}

\begin{rem}
It is well known that the orientation reversing homeomorphism of
$S^3$ is unique up to isotopy, but this is not true for the solid
tours. The ambiguity can be removed by posing the framing fixing
condition used in the above definition.
\end{rem}

To prove the main theorem of this subsection, we need the following
two lemmas.

\begin{lem}\label{lem:closed braid}
Suppose $l'$ is a framing of a thick braid in $N$. If $l'$ is
strictly achiral with respect to a framing $l$ of $N$, then the
strict achirality can be given by a $D^2$-fiberation preserving
homeomorphism of $N$.
\end{lem}

\begin{proof}
Suppose the strict achirality of $l'\subset N$ with respect to $l$
is defined by an orientation reversing homeomorphism $r : N \to N$.
Then we can isotope $r$ relative to $l' \cup l$ to a
$D^2$-fiberation preserving homeomorphism $r'$: since we can slide
$r(D^2 \times
*)$ relative to $l' \cup l$ to $D^2 \times *$, and the isotopy is
just the process of this sliding. (More clear way to see this
sliding: let $p:\tilde N\to N$ be the infinite cyclic covering,
and $\tilde l'\subset \tilde N$ be the preimage of $l'$, then
$(\tilde N, \tilde l')$ is homeomorphic to $(D^2, \text{$w$
points}) \times \R$ where $w$ is the winding number of $l'$ in
$N$.)
\end{proof}

\begin{lem}\label{lem:windding number}
Suppose $\Sigma\subset S^3$ is a tame embedding given by a
defining sequence $\{N_n\}_{n\ge0}$. If $\Sigma\subset \Fix(r)$
for some homeomorphism $r: S^3 \to S^3$, then there exists $k>0$
such that $r(N_n)\subset \intr N_0$, and moreover $r(N_n)$, $N_n$
have the same winding number in $N_0$ for $n\ge k$.
\end{lem}

\begin{proof} Since $S^3$ is compact, $r$ is uniformly continuous.

(i) Let $\epsilon = d(N_1, \partial{N}_0)/2$.

(ii) Choose $0<\delta<\epsilon$ such that if $d(x,x')<\delta$ then
$d(r(x),r(x'))<\epsilon$.

(iii) Choose $k>0$ such that $\max_{x \in N_k} d(x,\Sigma)<\delta$.

Now fix an integer $n\ge k$. For any $x\in N_n$, by (i) we have
$d(x,\partial{N}_0) \ge 2\epsilon$ and by (iii) we can choose
$x'\in\Sigma$ such that $d(x,x')<\delta$, hence by (ii) we have
$$d(x,r(x)) \le d(x,x') + d(x',r(x')) + d(r(x'),r(x))
  < \delta + 0 + \epsilon < 2\epsilon \le d(x,\partial{N}_0).
$$
It follows that the unique geodesic $\alpha(x)$ connecting $x$ and
$r(x)$ lies in $\intr N_0$. Therefore, $\{\alpha(x) \mid x \in
N_n\}$ gives rise to a homotopy from $N_n$ to $r(N_n)$ in $N_0$. In
particular, $r(N_n)\subset \intr N_0$ and $r(N_n),N_n$ have the same
winding number in $N_0$.
\end{proof}

\begin{thm}\label{thm:strictly achiral}
Let $\Sigma \subset S^3$ be a tame solenoid with defining sequence
$\{N_n\}_{n\ge0}$ and let $l_n$ denote a zero framing of $N_n$ in
$S^3$, that is, $l_n$ is null-homologous in $S^3\setminus N_n$. Then
$\Sigma$ is strictly achiral if and only if there exists $k\ge0$
such that the $l_k$ is strictly achiral in $S^3$, and $l_{n+1}$ is
strictly achiral in $N_n$ with respect to $l_n$ for all $n\ge k$.
\end{thm}

\begin{proof}
Sufficiency. Without loss of generality, we assume $k=0$. Since
$l_0$ is strictly achiral, there is an orientation reversing
homeomorphism $f_0 : (S^3,N_0) \to (S^3,N_0)$ such that $f_0|_{N_0}$
is $D^2$-fiberation preserving and fixes $l_0$ pointwisely.

Since $l_1$ is strictly achiral in $N_0$ with respect to $l_0$, by
Lemma \ref{lem:closed braid} there is an orientation reversing and
$D^2$-fiberation preserving homeomorphism $f_1 : (N_0,N_1) \to
(N_0,N_1)$ such that $l_1 \cup l_0$ stays in $\Fix(f_1)$. Since both
$f_0|_{N_0}$ and $f_1$ are orientation reversing and
$D^2$-fiberation preserving homeomorphisms of $N_0$ and both fix
$l_0$ pointwisely, we may assume $f_1|_{\partial{N}_0} =
f_0|_{\partial{N}_0}$. So $f_1$ may be extended onto $S^3$ by
setting $f_1|_{S^3\setminus N_0} = f_0|_{S^3\setminus N_0}$.

Then for $n>1$, by the same reason we can recursively define
homeomorphism $f_n : (S^3,N_n) \to (S^3,N_n)$ such that
$f_n|_{S^3\setminus N_{n-1}} = f_{n-1}|_{S^3\setminus N_{n-1}}$ and
that $f_n|_{N_{n-1}}$ is $D^2$-fiberation preserving and fixes $l_n
\cup l_{n-1}$ pointwisely.

To apply Lemma \ref{lem:converge}, set $U_n=N_n$. Clearly Conditions
(1) and (3) of Lemma \ref{lem:converge} are satisfied. Since
$f_n|_{N_n}$ is $D^2$-fiberation preserving and the diameters of the
meridian disks of $N_n$ tend to zero uniformly as $n\to\infty$,
Conditions (2) and (4) of Lemma \ref{lem:converge} are also
satisfied. By Lemma \ref{lem:converge}, $f_n$ uniformly converges to
an orientation reversing homeomorphism $f : S^3 \to S^3$ with the
property that $\cap_{n\ge0} N_n \subset \Fix(f)$.

Necessity. Suppose the strictly achirality of $\Sigma$ is defined by
an orientation reversing homeomorphism $r$ and let $k>0$ be given by
Lemma \ref{lem:windding number}. Then for any $n\ge k$ we have that
both $N_n$ and $r(N_n)$ are contained in the interior of $N_0$.

Fix an integer $n\ge k$ and choose a big integer $j$ so that $N_j
\subset N_n \cap r(N_n)$. By Lemma \ref{lem:braid} (3), both
$\partial{N}_n$, $r(\partial{N}_n)$ are isotopic in $N_0\setminus
N_j$ to some components of $\Gamma \cup \partial{N}_0 \cup
\partial{N}_j$ where $\Gamma$ is the JSJ-decomposition tori of
$N_0\setminus N_j$. Therefore, we can isotope $r$ with support in
$N_0\setminus\Sigma$ so that either $r(N_n)=N_n$, or by Lemma
\ref{lem:braid} (1) $r(N_n) \subset N_n$ or $N_n \subset r(N_n)$ is
a thick braid of winding number greater than 1. By Lemma
\ref{lem:windding number} $r(N_n),N_n$ have the same winding number
in $N_0$, so the latter case cannot happen and, moreover, we may
assume the zero framing $l_n$ of $N_n$ lies in $\Fix(r)$.

By the same argument, $r$ can be further isotoped with support in
$N_n\setminus\Sigma$ so that $r(N_{n+1})=N_{n+1}$ and $l_{n+1}
\subset \Fix(r)$. Therefore, $l_n$ is strictly achiral in $S^3$ and
$l_{n+1}$ is strictly achiral in $N_n$ with respect to $l_n$.
\end{proof}

\subsection{Examples}

Thanks to Theorem \ref{thm:strictly achiral}, the strict achirality
of tame solenoids is the problem of strict achirality of knots in
$S^3$ and closed braids in the solid torus. Below we fix a point $*
\in \partial{D}^2$ and let $l_*$ denote the framing $*\times S^1$ of
$D^2\times S^1$.

\begin{defn}
Call a braid $\beta$ {\it achiral}, if $\beta$ is conjugate to its
mirror image $\beta^*$.
\end{defn}

Note that the achirality is well defined for a conjugacy class of
braids. Also note that the closure $\overline\beta \subset D^2\times
S^1$ of a braid $\beta$ is connected if and only if $\beta$ is
cyclic, i.e. $\beta$ permutes its strands cyclicly.

\begin{prop}\label{prop:achiral braid}
For every cyclic braid $\beta$, the closure $\overline\beta\subset
D^2\times S^1$ is strictly achiral with respect to $l_*$ if and only
if $\beta$ is achiral.
\end{prop}

\begin{proof}
Sufficiency. Suppose $\beta = \alpha^{-1} \beta^* \alpha$. Then
$\beta = (\alpha\beta^k)^{-1} \beta^* (\alpha\beta^k)$ for any
integer $k$. Since $\beta$ and $\beta^*$ give rise to the same
cyclic permutation on their strands, it follows from the equality
$\beta = \alpha^{-1} \beta^* \alpha$ that the permutation given by
$\alpha$ is commutative with hence is a power of that given by
$\beta$. So we can choose $k$ so that $\alpha\beta^k$ is a pure
braid. Replacing $\alpha$ by $\alpha\beta^k$, we may assume $\alpha$
is a pure braid.

By the definition of braid, there is a $D^2$-fiberation preserving
and boundary fixing homeomorphism $\tilde g: D^2\times [0,1]\to
D^2\times [0,1]$ which sends $\alpha^{-1} \beta^* \alpha$ to
$\beta$. Then $\tilde g$ induces a homeomorphism $g: D^2\times
S^1\to D^2\times S^1$ which sends $\overline {\alpha^{-1} \beta^*
\alpha}$ to $\overline\beta$. Moreover, by canceling the $\alpha$
part and the $\alpha^{-1}$ part of $\overline{\alpha^{-1} \beta^*
\alpha}$, we can define a $D^2$-fiberation preserving and boundary
fixing homeomorphism $f: D^2\times S^1\to D^2\times S^1$ which sends
$\overline{\beta^*}$ to $\overline {\alpha^{-1} \beta^* \alpha}$.

Finally, let $r: D^2\times S^1\to D^2\times S^1$ be the reflection
about the ``page'' (assume $l_*$ lies on it). Then $gfr: D^2\times
S^1\to D^2\times S^1$ is an orientation reversing homeomorphism and
sends $\overline\beta$ to $\overline\beta$. Since $\alpha$ is a pure
braid and all $g,f,r$ preserve each $D^2$ slice of $D^2\times S^1$,
it follows that $gfr$ fixes both the closure $\overline\beta$ and
the framing $l_*$ pointwisely.

\bigskip
\begin{center}
    \psfrag{r}[]{$r$}
    \psfrag{f}[]{$f$}
    \psfrag{g}[]{$g$}
    \includegraphics[scale=.8]{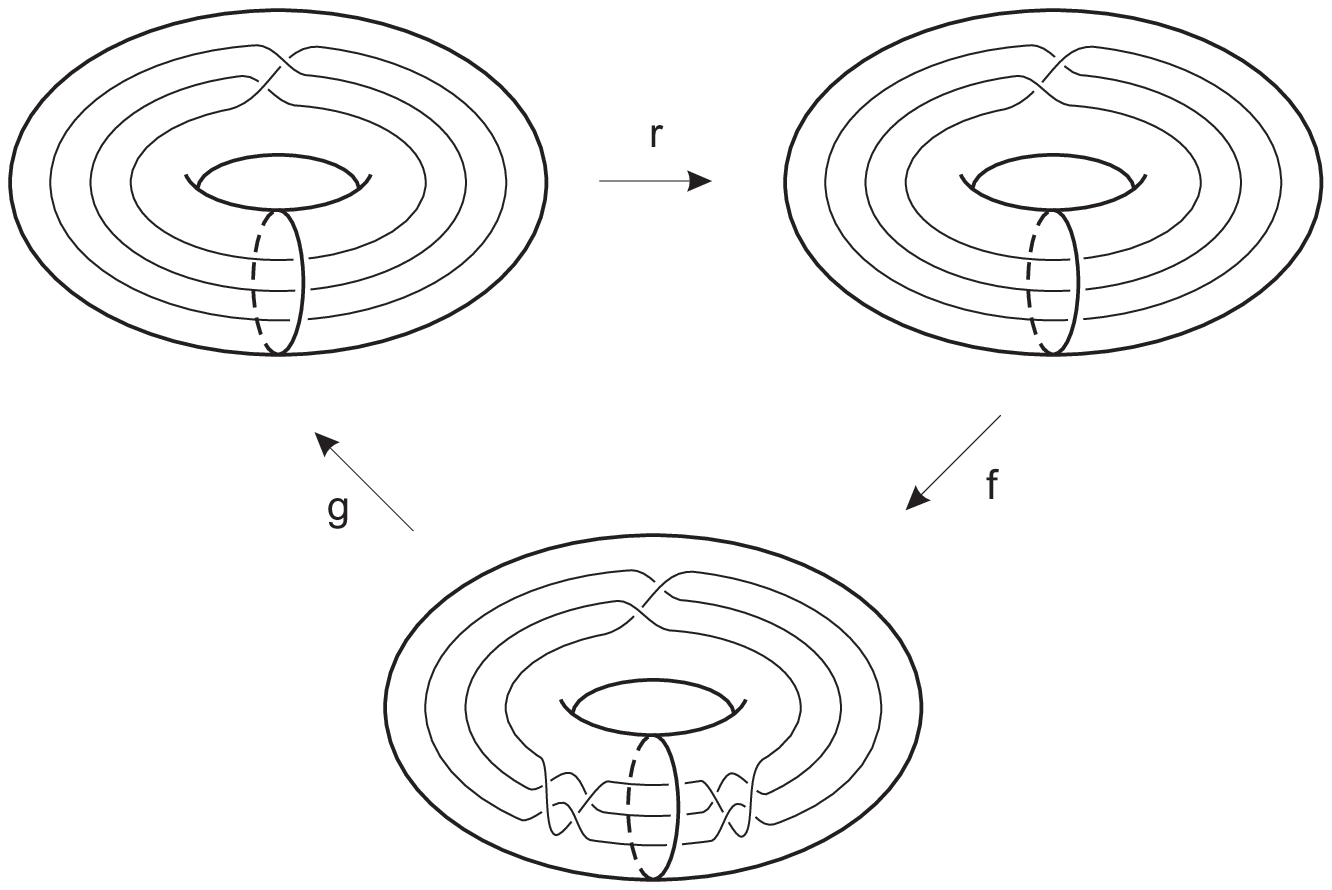}
\end{center}
\centerline{Figure 4}
\bigskip

Figure 4 illustrates the case $\beta=\sigma_1\sigma_2^{-1}$ and
$\alpha=\sigma_2\sigma_1^2\sigma_2^{-1}$ (see Example
\ref{exam:achrial braid} for the equality $\beta=\alpha^{-1} \beta^*
\alpha$).

Necessity. If $\overline\beta$ is strictly achiral with respect to
$l_*$ then $\overline\beta$ is isotopic to its mirror image,
therefore $\beta$ is conjugate to its mirror image $\beta^*$.
\end{proof}

\begin{exam}\label{exam:achrial braid}
Examples of cyclic, achiral braids.

(1) $\beta=\sigma_1\sigma_2^{-1}$ is cyclic and achiral. Setting
$\alpha = \sigma_2\sigma_1^2\sigma_2^{-1}$, one can verify the
equality $\beta=\alpha^{-1}\beta^*\alpha$ either by directly a braid
move or by the substitution of braid relation
$\sigma_2\sigma_1\sigma_2=\sigma_1\sigma_2\sigma_1$ as follows
\begin{align*}
  \alpha^{-1}\beta^*\alpha
  & = (\sigma_2\sigma_1^{-1}\sigma_1^{-1}\sigma_2^{-1}) (\sigma_1^{-1}\sigma_2)
  (\sigma_2\sigma_1\sigma_1\sigma_2^{-1}) \\
  & = \sigma_2\sigma_1^{-1}\sigma_2^{-1}\sigma_1^{-1} \sigma_2^{-1}\sigma_2
  \sigma_2\sigma_1\sigma_1\sigma_2^{-1} \\
  & = \sigma_2\sigma_2^{-1}\sigma_1^{-1}\sigma_2^{-1} \sigma_2^{-1}\sigma_2
  \sigma_2\sigma_1\sigma_1\sigma_2^{-1} \\
  & = \sigma_1\sigma_2^{-1} = \beta.
\end{align*}

(2) For any braid $\beta$, $\beta \beta^*$ is achiral, since
$${\beta^*}^{-1}(\beta \beta^*)^*\beta^* = {\beta^*}^{-1}\beta^* \beta\beta^* = \beta\beta^*.$$
Moreover, for any cyclic braid $\beta$ of odd number of strands,
$\beta\beta^*$ is also cyclic. Hence for each cyclic braid $\beta$
of odd number of strands, $\beta\beta^*$ is cyclic and achiral.

(3) If $\beta$ is an achiral braid, then so is $\beta^k$ for any
integer $k$. Moreover, if $\beta$ is cyclic, then so is $\beta^k$
for every integer $k$ relatively prime to the number of strands.
\end{exam}

\begin{prop}\label{prop:chiral braid}
(1) If a connected closed braid $\overline\beta\subset D^2\times
S^1$ is strictly achiral with respect to $l_*$, then the writhe of
$\overline\beta$ is zero (for the definition of writhe, see [A,
p.152]).

(2) Hence a connected closed braid $\overline\beta\subset D^2\times
S^1$ is not strictly achiral with respect to $l_*$ if
$\overline\beta$ is either of even winding number, or a cable.
\end{prop}

\begin{proof}
(1) By Proposition \ref{prop:achiral braid}, $\overline\beta$ is
strictly achiral implies $\beta=\alpha^{-1}\beta^*\alpha$ for some
braid $\alpha$. Clearly $wr(\overline\beta) =
wr(\overline{\alpha^{-1}\beta^*\alpha}) = wr(\overline{\beta^*}) =
-wr(\overline\beta)$. It follows that $wr(\overline\beta)=0$.

(2) Suppose $\overline\beta$ is connected and is of even winding
number. It is an elementary exercise to show that the number of
the crossings of $\overline\beta$ is odd. Hence
$wr(\overline\beta)$ must be odd, which contradicts (1).

Suppose $\overline\beta$ is a cable. Then all the crossings of
$\overline\beta$ have the same sign, so $wr(\overline\beta)$ is
non-zero, which contradicts (1).
\end{proof}

Now we state the main result of this subsection.

\begin{thm}\label{thm:sa embedding}
A solenoid of type $\varpi=(w_1,w_2,...,w_n,...)$ has a strictly
achiral tame embedding into $S^3$ if and only if all except finitely
many $w_n$ are odd.
\end{thm}

\begin{proof}
Necessity is immediate from Theorem \ref{thm:strictly achiral} and
Proposition \ref{prop:chiral braid} (2).

Sufficiency. Assume all $w_n$ are odd. Let $N_0$ be a tubular
neighborhood of a strictly achiral knot in $S^3$ and let $N_n$ be a
tubular neighborhood of $\overline{\beta_n\beta_n^*}$ in $N_{n-1}$
where $\beta_n$ is an arbitrary cyclic braid on $w_n$ strands. By
Theorem \ref{thm:strictly achiral}, Proposition \ref{prop:achiral
braid} and Example \ref{exam:achrial braid} (2) the defining
sequence $\{N_n\}_{n\ge0}$ gives rise to a strictly achiral tame
embedding of the solenoid of type $\varpi=(w_1,w_2,...,w_n,...)$.
\end{proof}

\begin{exam}
(1) The 2-adic solenoid has no strictly achiral tame embedding into
$S^3$.

(2) By Lemma \ref{lem:MP} (2), Example \ref{exam:achrial braid} (1)
and Proposition \ref{prop:chiral braid} (2), up to equivalence the
3-adic solenoid has a unique unknotted strictly achiral tame
embedding into $S^3$, which is yielded by nesting the thick braid
pictured in the middle of Figure 3.

(3) The embedding in Theorem \ref{thm:sa embedding} can be chosen to
be either knotted or unknotted, by letting in the proof either $N_0$
be a tubular neighborhood of the figure-8 knot, or letting $N_0$ be
a tubular neighborhood of the unknot and $\beta_n =
\sigma_1\sigma_2\cdots\sigma_{w_n-1}$ (we leave it to the reader to
verify that the closure $\overline{\beta_n\beta_n^*}$ is unknotted
in $S^3$).
\end{exam}

\bibliographystyle{amsalpha}

\end{document}